# Some Ideas for Improving Stock Price Prediction Based on Machine Learning

Negin Bagherpour

negin.bagherpour@ut.ac.ir

Assistant Professor, University of Tehran, Department of Engineering Sciences

**Abstract:** Stock price prediction is a complicated and interesting task. Noisy trends make stock pricing sensitive and complicated while the economical motivation behind, keeps it interesting for researchers and investors. In this paper we are to outline two novel ideas for stock pricing. We also test each of our suggested algorithms for predicting the price of 6 stocks from different sectors. To show the efficiency of our proposed algorithm, we compare the predicted prices with real values and also perform a backtest to verify that the annual returns based on real data and predicted price are almost the same.

1. Introduction

Stock price prediction is very difficult due to the fact that the price trend is usually very noisy and is affected by too many factors including economic growth rate, currency market and unemployment rate [1] Two possible strategies come into mind to compute a reliable prediction: 1) Suggesting a complicated nonlinear model to make a precise prediction and 2) Improving the predicted price by simple models such as linear regression. In this paper the second strategy has been considered to suggest a reliable stock price prediction.

Predicting market price is useful for investors. The tools for the stock market prediction can monitor, predict and regulate the market to guide the investors to correct decisions. The methods for price estimation are mainly divided into two types: prediction and classification. In prediction basically a model is suggested to estimate the price of each stock in future. The sufficiency might be certified based on the past price data. In classification however a prediction is computed for the price trend and stocks with similar sale and purchase signals are located in same cluster. The performance is examined based on average cluster accuracy. These tools help the investors to predict future behaviors and make reliable sale/purchase decisions based on intelligent data analysis methods such as regression approach [2], Knowledge Discovery in Databases (KDD) [3] and fuzzy models [4]. These data mining techniques are basically essential for searching the hidden parts and increase the accuracy level. The prediction methods may be divided into econometric and machine learning models. In 1994, Booth et. al. [5] used six explanatory variables to outline the ARIMA econometric model to predict stock price. In 1998 GARCH model was introduced [6]. As a long memory stochastic volatility model, it was superior to ARIMA and other existing models. Zhang et al. proposed a hybrid ARIMA-GARCH model based on differential information in 2018 [7]. Adding the approximate differential information of the dependent variable lag and considering stock price change trend information, improved the price trend prediction. Development of machine learning motivated the researchers trying to present more accurate predictions by use of this new technology instead of traditional prediction models. Maknickas and Maknickiene [8] used a recursive neural network (RNN) to provide a stock price prediction model and analyzed the optimal values for RNN parameters, such as the number of neurons and the number of iterations. The LSTM model was also used to predict the rise and fall of the stock price in 2017 [9]. LSTM showed to suggest more accurate price estimations than that of the traditional machine learning model and non-time series model. Moreover, in 2019 Peng et al. [10] focused on the preprocessing methods such as interpolation, wavelet denoising and normalization of data and tried various parameter values in LSTM model. They found that the optimized model significantly improved the prediction accuracy while the computational complexity was low. In addition, Vo et al. [11] studied Bi-directional LSTM (Bi-LSTM) and found that it read the data one more time leading to prediction accuracy improvement. On the other hand, a combination of statistical econometric models with machine learning models or use of more than two different machine learning models at the same time to predict stock price is also common. These models usually have better performance than single ones. Achkar et al. [12] considered two different model combination methods including back propagation-multi layer perception (BPA-MLP) and LSTM-RNN in 2018. For the price data of Facebook, Google and Bitcoin, LSTM-RNN model worked better than BPA-MLP [12]. Bao et al. [13] suggested a two-step mechanism: noise reduction in original time series stock data through wavelet transform and price prediction by LSTM model. M'ng and Mehralizadeh [14] proposed a prediction model named wavelet principal component analysis-neural network (WPCA-NN), which combined wavelet transform, PCA and artificial neural network to denoising, removing the random noise in the stock price sequence. A hybrid of PCA and LSTM was suggested in [15]. The results showed that the performance of the WPCA-LSTM was better than traditional prediction methods. In 2019 an LSTM-CNN [16] model based on feature combination, using stock time series and stock trend graphs as input features was also created. Based on the reported numerical results the LSTM-CNN model was superior to the single model in predicting stock prices. On the other hand, DM (distributed memory) and DBOW (distributed bag of words) methods were suggested and a successful combination was outlined in [17]. This hybrid was also tested besides stacked auto-encoder (SAE) [18]. Due to fluctuating trends in stock pricing last available data might more effectively predict future price. More precisely, the price in the last few days would have more influence on the future price and it looks reasonable to update the model based on the last days. Online learning algorithms can address these problems with incremental learning approaches. Recently, online learning algorithms was used abundantly for stock pricing; see e.g. [19], [20]. A bothering problem in most of the online machine learning techniques was the accumulation of unnecessary part of data in the past. In fact, feature selection is not possible in real time [21]. The ARTML (Adaptive Real-Time Machine Learning) library [22] was developed in python with the following flexibilities:

- Incremental learning - Updating the model with new data in real time
- Decremental learning - Updating the model by excluding unnecessary data in real time
- Real time attribute addition - Adding new features to the model in real time
- Real time attribute deletion - Deleting features from model in real time
- Distribute processing - Power to build models using distributed data
- Parallel processing - carrying out parallel processing using GPU's

In summary, there exists a variety of methods for stock price prediction including simple and complicated algorithms. The complicated methods are based on econometric or machine learning techniques. Although more complicated methods such as LSTM-RNN and LSTM-CNN might suggest an accurate prediction, high computational complexity is inevitable. Hence, there is a motivation behind developing new algorithms which are precise enough without needing a lot of computations. These algorithms are known as simple methods. In this paper we develop two algorithms to compute a price prediction based on incremental learning which is more accurate than any set of existing simple method with low computational complexity. By use of these two algorithms, a prediction which is more accurate than simple methods such as linear regression and moving average is in hand. We note that these algorithms might also be used to improve more complicated algorithms but it is not efficient because of high computational complexity of initial methods and lower possible improvement.

The rest of this paper is organized as follows: in Section 2 we present a brief review on well-known machine learning algorithms for stock price prediction. In sections 3 and 4 two new algorithms are outlined for presenting a more reliable pricing. Numerical tests for pricing 6 stocks belonging to different sectors is presented in Section 5. Finally, Section 6 gives our concluding remarks.

2. Machine Learning Techniques for Stock Pricing

Stock pricing is known as an interesting and challenging context in recent years. It is interesting because of economical motivation behind price prediction and it is challenging because of nonlinear and very noisy behavior of values. There exists a variety of methods for stock pricing many of which are in machine learning. Here, we briefly review the main simple methods in machine learning [23]. (The complicated methods were reviewed in Section 1 and are not directly used in the rest of this research):

(a) Linear Regression (LR): The best linear model based on least squares formulation for each $k$ days is considered to predict the price of next day. For a normal price trend $k$ equal to 4 or 5 seems reasonable. LR works well as long as we have no sudden fluctuation. However, any intense change makes the locally linear regression inaccurate.

(b) Moving Average (MA) [24]: The average price of each $k$ days is assumed as an approximation of the next day price. MA also works well for smooth trends and for any sudden change, the average price of a few days may be considerably far from the next day price.

(c) Support Vector Regression (SVR) [25]: Support vector Machine (SVM) is a robust prediction method for data classification based on supervised learning models. Given a set of training examples belonging two different categories, an SVM training algorithm builds a model that assigns new examples to one the two groups. The main idea of the SVM algorithm is to map each training example to a point in space by use of kernel functions, so as to maximize the width of the gap between the two categories to keep them significantly separated. New examples would be mapped into that same space and predicted to belong to one category based on the side of the gap they fall. Support vector regression (SVR) is characterized by the use of kernels, sparse solution and control of the margin and the number of support vectors to fit a proper regression model. Although SVR is not as common as SVM, it is an effective tool for approximating real-valued functions.

(d) Pair Trading (PT): The dependency of each two stock prices is reviewed by use of statistical tests and some beneficial trades are suggested between each to dependent stocks. The most important property of PT is its ability to present a trade strategy with proper return without encountering complicated computations.

(e) Nonlinear Regression (NR): A nonlinear function of different factors is fitted to estimate each day price. Normally, stock price is assumed to be a function of economic growth rate, currency value, gold price, average income and unemployment rate. Based on the history of price, the dependency type (periodic, polynomial, exponential, logarithmic and …) to each factor would be revealed. The coefficients of the suggested terms might be computed by solving the corresponding nonlinear least squares problem.

(f) Incremental Lear

In sections 3 and 4 two new strategies are discussed to suggest a more accurate prediction (with lower average error) than a group of available methods. These strategies improve existing predictions by very low computational cost.

3. Error Analysis Strategy

Here, we aim to outline a new algorithm for suggesting the best simple prediction method in each day. We first note that a complicated nonlinear model, despite the higher computational complexity, doesn't necessarily provide an accurate market prediction. In fact, a more efficient idea is to suggest a procedure to indicate the best estimation among any set of simple prediction methods such as MA, LR and SVR. To do this, we define success ratio for each prediction method $M$, in a period of $k$ days SR$(M, k)$, to be the number of days for which the method $M$ leads to the best estimation. In each day, the method with the maximum success ratio in the last $k$ days is chosen to predict the next day price but it is not efficient to compute the success ratio daily. Hence, we compute the success ratio every $p$ days and start price prediction using the method with greatest success ratio $M$. When the success ratio of $M$ decreases, it is computed in a shorter period to see whether the reduction was due to a serious trend change. In case the success ratio in the shorter period is not lower than the previous value, price prediction continues based on $M$. Otherwise, the method $N$ with the second largest success ratio would be selected to continue estimation. The motivation behind this idea is that reduction of success ratio might occur because of a sudden and temporary trend change or because of a gradual change. We note that in case of a gradual change it is necessary to change the prediction method while in a sudden change it is not. The outlined procedure is actually investigation of derivative changes of success ratio. The main steps of our proposed algorithm are summarized in Figure 1.

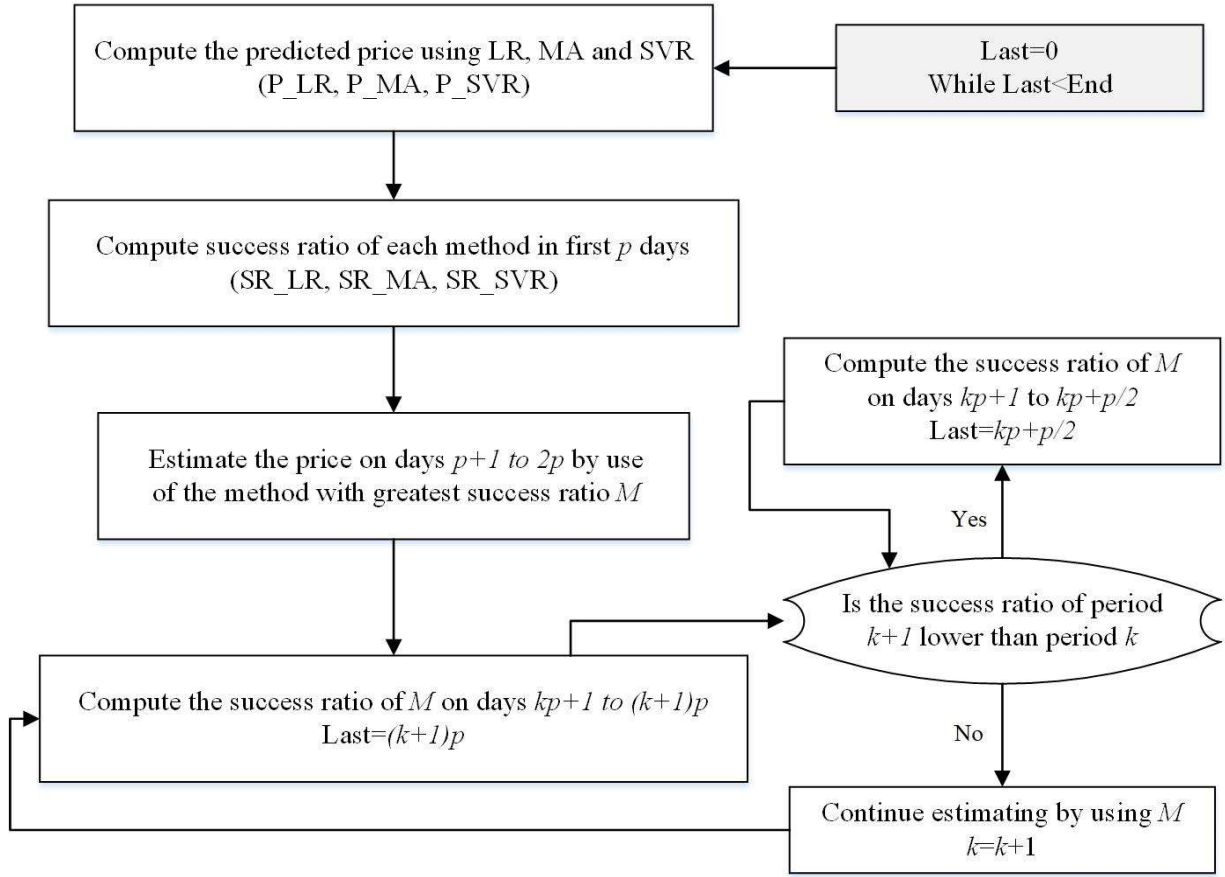

Figure 1: Error Analysis Strategy for Price Prediction

In noisy trends, error analysis leads to a sufficiently fine approximation by smoothing the fluctuations; however, the estimation for relatively smooth trends is not as accurate as expected. In the following we suggest a new strategy for relatively smooth data.

4. Error Regression Strategy

In this strategy, we first assign an index to each method and we interpolate the index of best prediction method in each day based on the last few days. We refer to this strategy as Error Regression (ER). Since the best prediction method in every single day affects the suggested method in ER, it is more accurate than error analysis. We note that to provide an easy to compute interpolant, discrete values might be chosen for each method; see Table 1. Let correspond an index $i$ to each of the $m$ methods. The resulting data for best method in each day would be as in Table 1.

Table 1: Data structure for best prediction method in each day (Example)

| Method | Method 1 | Method 2 | Method 3 | Method 4 |
|---|---|---|---|---|
| Index | 1 | 2 | 3 | 4 |
| Day 1 |   |   |   |   |
| Day 2 |   |   |   |   |
| Day 3 |   |   |   |   |

The data in Table 1 leads to the following points for interpolation:

$$d=1, I=1$$
$$d=2, I=1$$
$$d=3, I=4$$

and the interpolant polynomial would be $I(d) = \frac{3}{2}d^2 - \frac{9}{2}d + 4$ which for instance predict that the index of best method for forth day is 10. There are two computational problems:

1) The resulting index might not be an integer. In this case, we might substitute the integer part.

2) The resulting index might not be in the existing domain. To solve this issue,
   - we might use bounded interpolants such trigonometric functions or
   - we might add bound conditions.

5. Numerical Tests

In this section we provide numerical test results to confirm the efficiency of both suggested strategies error analysis (EA) and error regression (ER). To this end, we first report the results for normal stocks in Section 5.1. We then represent the results for stocks with smoother trends in Section 5.2. The predictions are compared to real prices based on two reasonable criteria:

1) Backtest to compare the resulting return,
2) Closeness of predicted and real prices.

Moreover, in Section 5.3 EA and ER are compared to suggest the user when and how much ER outperforms EA.

5.1. Error Analysis Algorithm

Here, six stocks including Apple, Goldman Sachs, Johnson & Johnson, Moderna, Nike and Pfizer are selected to be estimated in last year by use of DM-DBOW [17], DM-DBOW-SAE [18], DOC-W-LSTM [18], LSTM-F [18], PCA-LSTM [15] and EA. The results are plotted in Figure 2.

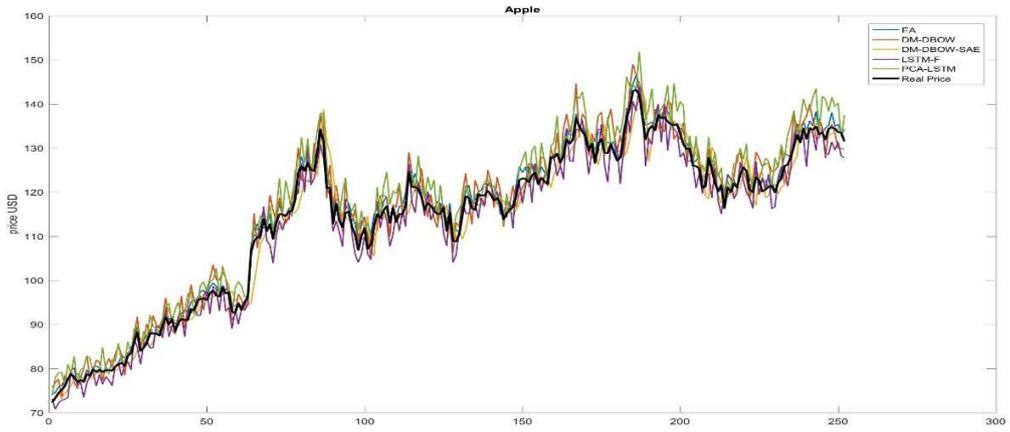

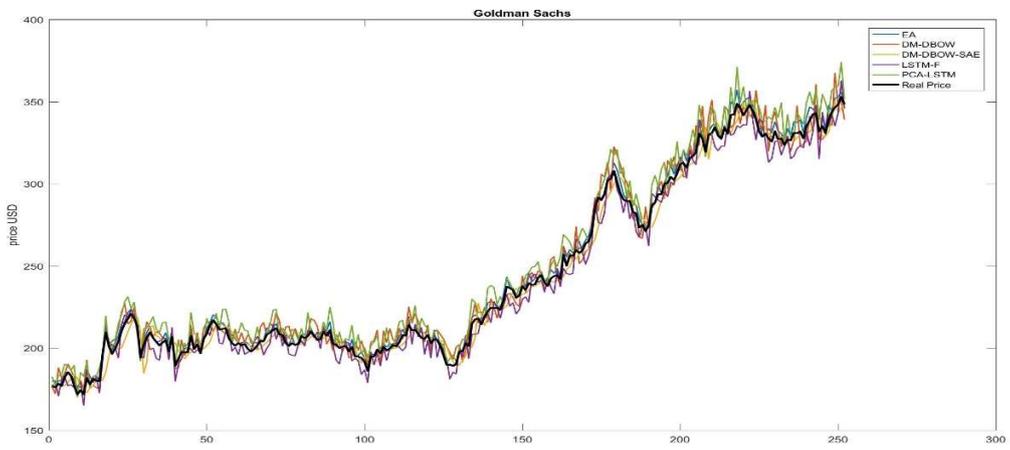

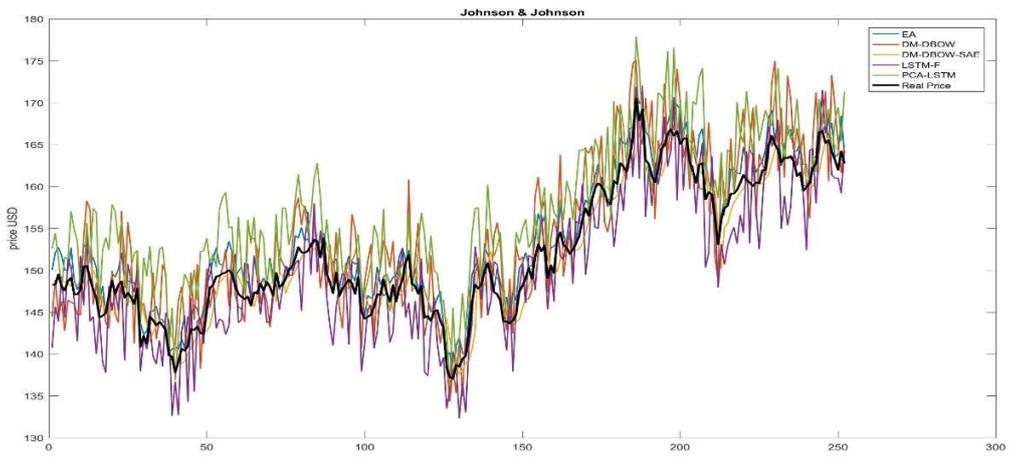

Figure 2: Price prediction for six selected stocks

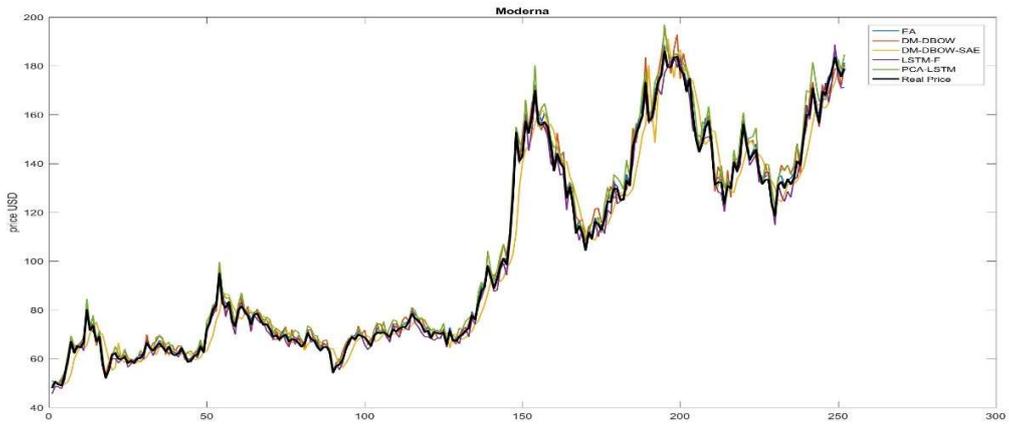

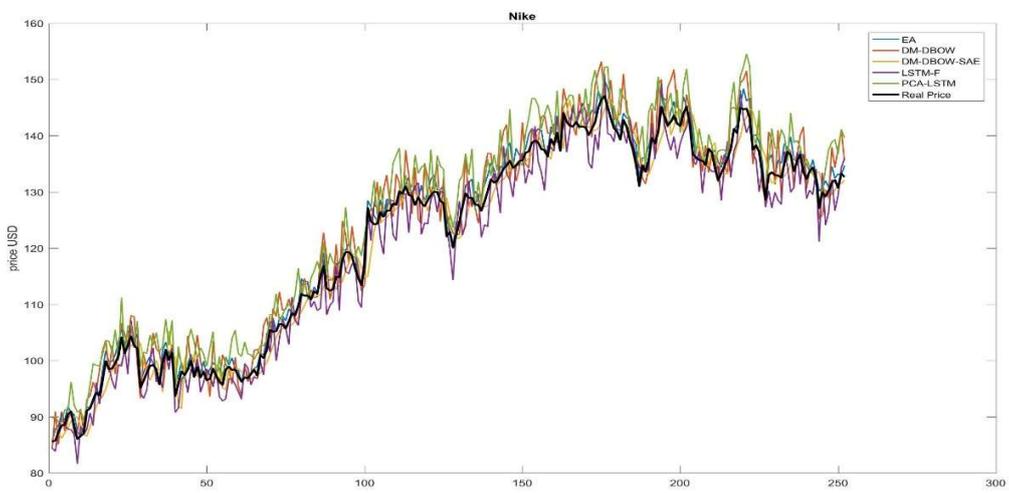

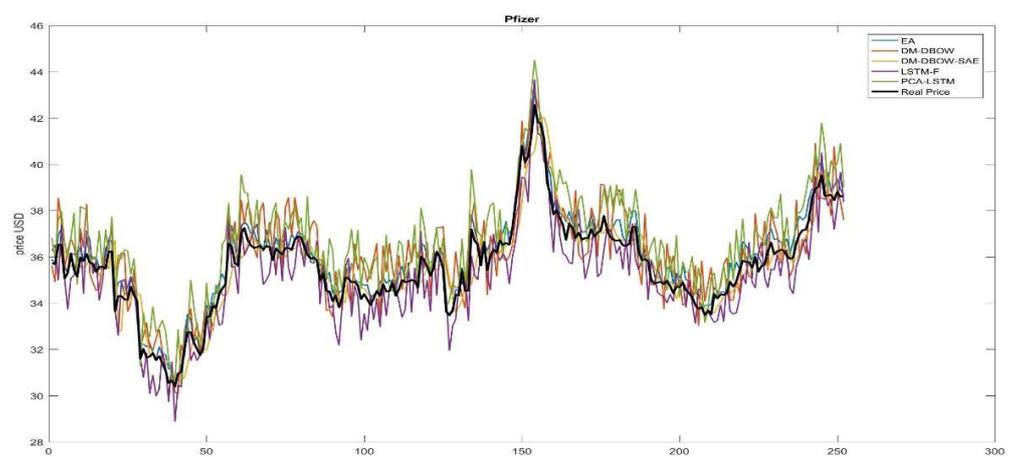

Figure 2: Price prediction for six selected stocks (Cont'd)

To confirm the effectiveness of EA strategy to improve MA and LR results, the average absolute and relative errors are reported in Table 2. Based on the average error values, the best prediction is provided by EA for all six stocks. Moreover, MA performed better than LR for Apple, Goldman Sachs, Johnson & Johnson and Moderna while LR outperforms MA for Nike and Pfizer.

Table 2: Average absolute/relative error values for price prediction based on EA, LR and MA

|  | EA(Abs/RMSE) | | LR(Abs/RMSE) | | MA(Abs/RMSE) | |
| --- | --- | --- | --- | --- | --- | --- |
| Apple | 2.2678 | 0.01725 | 2.5059 | 0.01906 | 2.3801 | 0.01811 |
| Goldman Sachs | 4.4568 | 0.01279 | 4.8977 | 0.01405 | 4.7249 | 0.01356 |
| Johnson & Johnson | 1.5013 | 0.00923 | 1.6954 | 0.01042 | 1.6088 | 0.00988 |
| Moderna | 4.9817 | 0.02786 | 5.2727 | 0.02949 | 5.2013 | 0.02909 |
| Nike | 1.9750 | 0.01489 | 2.0102 | 0.01516 | 2.0996 | 0.01583 |
| Pfizer | 0.4849 | 0.01255 | 0.5141 | 0.01331 | 0.5176 | 0.01339 |

As a second numerical test to confirm the efficiency of EA, a backtest is done to estimate the annual return values. As represented in Table 3, the annual returns computed based on real prices and the estimated prices by EA are almost the same.

Table 3: Annual Return % by Backtest on EA, MA and LR

|  | Real Price | EA | LR | MA |
| --- | --- | --- | --- | --- |
| Apple | 16.7 | 16.7 | 11.2 | 13.4 |
| Goldman Sachs | 15.3 | 15.2 | 14.1 | 14.8 |
| Johnson & Johnson | 18.6 | 18.6 | 17.4 | 18.1 |
| Moderna | 19.3 | 19.4 | 18.2 | 18.8 |
| Nike | 11.1 | 11.3 | 11.7 | 10.6 |
| Pfizer | 19.7 | 19.7 | 18.7 | 19.1 |

In Figure 2, it is shown that for real prices of Apple and predicted prices by EA the buy/sell signals are in accordance. This is a desirable feature of EA which can be easily established for all six stocks.

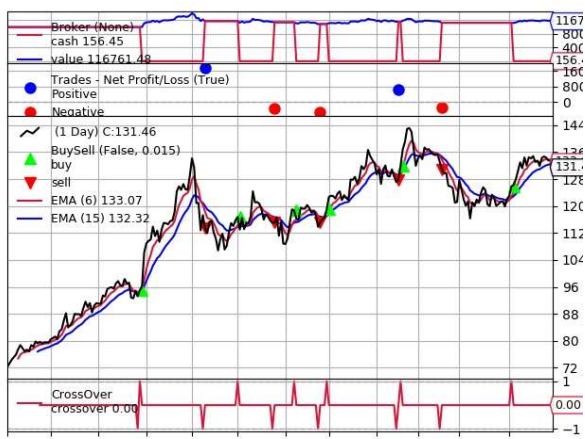
(a) Real Price

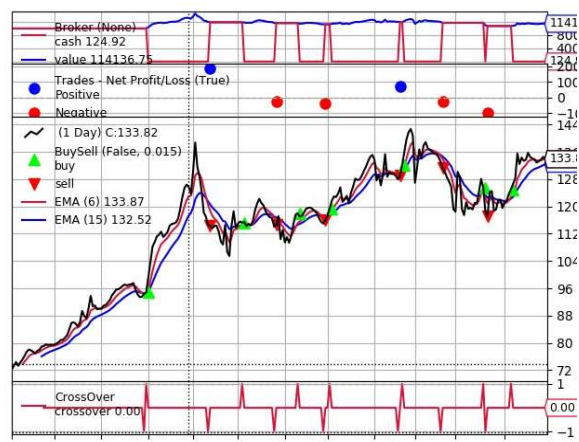
(b) EA

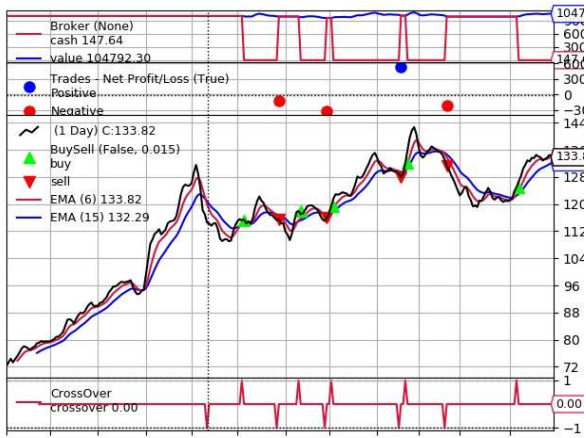
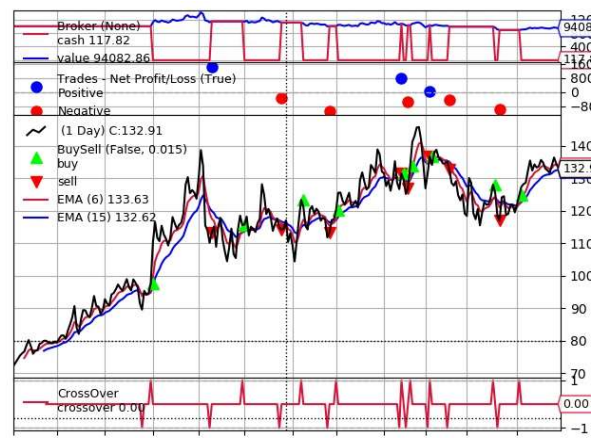

(c) MA
(d) LR

Figure 2: Backtest Strategies for Apple Real Price, EA, LR and MA

5.2. Error Regression Algorithm

In the following we provide numerical results of Error Regression (ER) algorithm for suggesting a proper prediction for the group of stocks with smoother trends. Hence, six stocks with relatively smooth trends including AbCellera Biologics, Cricut, DoorDash, Magnite, MediaAlpha and Singapore Airlines are chosen to check the efficiency of ER algorithm. In Table 4 the average absolute and relative error values are reported to make a more accurate comparison.

Table 4: Average absolute/relative error values for price prediction based on EA, LR and MA

|  | ER(Abs/Rel) | | LR(Abs/Rel) | | MA(Abs/Rel) | |
| --- | --- | --- | --- | --- | --- | --- |
| ABCL | 2.0311 | 0.06397 | 2.1420 | 0.06746 | 2.1744 | 0.06848 |
| Cricut | 0.7459 | 0.03039 | 0.9355 | 0.03812 | 0.9804 | 0.03995 |
| DoorDash | 7.2815 | 0.05162 | 7.7481 | 0.05492 | 7.3365 | 0.05201 |
| Magnite | 1.3242 | 0.05299 | 1.3575 | 0.05432 | 1.4938 | 0.05977 |
| MediaAlpha | 1.9671 | 0.05652 | 2.1406 | 0.06151 | 2.4378 | 0.07005 |
| Singapore Airlines | 0.0648 | 0.01873 | 0.0657 | 0.01899 | 0.0668 | 0.01931 |

Finally, we present the backtest results to compare the estimated annual return for the real prices, ER, LR and MA. Backtest results also confirm the efficiency of ER in suggesting similar trend to real data. The backtest results are represented in Table 5 and Figure 5. This is also of interest to check out whether EA and ER works significantly different for the selected stocks which is discussed in Section 5.3.

5.3. Comparison on EA and ER

Here, we present a comparison on EA and ER to identify their features. The average error values are reported in Table 6.

Table 6: Relative Error Values for EA and ER

| Group 1 | EA | ER | Group2 | EA | ER |
|---|---|---|---|---|---|
| **Apple** | 0.01725 | 0.01863 | **ABCL** | 0.07621 | 0.06397 |
| **Goldman Sachs** | 0.01279 | 0.01335 | **Cricut** | 0.03995 | 0.03039 |
| **Johnson & Johnson** | 0.00923 | 0.00978 | **DoorDash** | 0.05258 | 0.05162 |
| **Moderna** | 0.02786 | 0.03096 | **Magnite** | 0.06153 | 0.05299 |
| **Nike** | 0.01489 | 0.01486 | **MediaAlpha** | 0.07427 | 0.05652 |
| **Pfizer** | 0.01255 | 0.01247 | **Singapore Airlines** | 0.02068 | 0.01873 |

It can be seen in Table 6 that for stocks in Group 1, EA outperforms ER while in Group 2, ER looks better than EA as was expected analytically.

6. Concluding Remarks

Stock pricing is of much interest because of its applications in investing. In recent years Machine Learning shows to be one of the most powerful tools for predicting price. In this paper we outlined two new strategies for improving the results achieved by two well-known algorithms in machine learning, linear regression and moving average. The first idea is to select the most successful method in the last few days for predicting next day price. The second strategy is to interpolate the index of best method in the last few days and use the interpolant to compute the index best possible method to predict next day. These two ideas are named Error Analysis and Error Regression respectively. We also present a variety of numerical tests to confirm the efficiency of presented ideas in suggesting a more accurate estimation of price for 12 selected stocks. Moreover, Error Regression shows to be more precise for smoother trends because the best method in each single day comes into account.